\documentclass[11pt,twoside]{article}
\usepackage{latexsym}
\usepackage{amssymb,amsbsy,amsmath,amsfonts,amssymb,amscd}
\usepackage{graphicx}
\usepackage{hyperref}
\usepackage{color,yfonts}
\newcommand{\commentout}[1]{}
\newcommand{\R}{\mathbb{R}}

\newcommand {\al} {\alpha}
\newcommand {\e}  {\varepsilon}

\newcommand {\sg} {\sigma}

\newcommand {\Chi} {{\bf \raise 2pt \hbox{$\chi$}} }
\newcommand {\f}   {\frac}
\newcommand {\p}   {\partial}
\newcommand{\dis}{\displaystyle}

\newcommand{\beq}{\begin{equation}}
\newcommand{\eeq}{\end{equation}}
\newcommand{\bea} {\begin{array}{rl}}
\newcommand{\eea} {\end{array}}
\newcommand{\bepa}{\left\{ \begin{array}{l}}
\newcommand{\eepa} {\end{array}\right.}
\newtheorem{theorem}{Theorem}[section]

\newtheorem{definition}[theorem]{Definition}
\newtheorem{remark}[theorem]{Remark}

\newtheorem{proposition}[theorem]{Proposition}

\newcommand{\qed}{{ \hfill
                       {\unskip\kern 6pt\penalty 500 \raise -2pt\hbox{\vrule\vbox to 6pt{\hrule width 6pt
                       \vfill\hrule}\vrule} \par}   }}
\title{\Large \bf Transversal instability for the thermodiffusive reaction-diffusion system}

\author{
Michal Kowalczyk\thanks{Departamento de Ingenier'a Matem‡tica and Centro de Modelamiento Matem\'atico (UMI CNRS 2807), Universidad de Chile, Santiago, Chile. 
Email~:~kowalczy@dim.uchile.cl. M.K. and B.P.~have been partially supported by the FONDECYT grant 1130126, the ECOS project \emph{C11E07} and \emph{Fondo Basal CMM}}
\and Beno\^ \i t Perthame\thanks{Sorbonne Universit\'es, UPMC Univ Paris 06, UMR 7598, Laboratoire Jacques-Louis Lions, F-75005, Paris, France. Emails~:~benoit.perthame@upmc.fr,~nicolas.vauchelet@upmc.fr. B.P. and N.V. are supported by the french "ANR blanche" project Kibord:  ANR-13-BS01-0004} 
\thanks{CNRS, UMR 7598, Laboratoire Jacques-Louis Lions, F-75005, Paris, France}
\thanks{INRIA-Paris-Rocquencourt, EPC MAMBA, Domaine de Voluceau, BP105, 78153 Le Chesnay Cedex, France}
\and Nicolas Vauchelet\footnotemark[2] \footnotemark[3] \footnotemark[4]
}

\date{\today}

\begin{document}
\maketitle
\pagestyle{plain}
%\tableofcontents
\pagenumbering{arabic}

\begin{abstract}
The propagation of unstable interfaces is at the origin of remarkable patterns that are observed in various areas of science as chemical reactions, phase transitions, growth of bacterial colonies. Since a scalar equation generates usually stable waves, the simplest mathematical description relies on two by two reaction-diffusion systems. Our interest is the extension of the Fisher/KPP equation to a  two species reaction which represents  reactant concentration and temperature when used for flame propagation, bacterial population and nutrient concentration when used in biology. 

We study in which circumstances instabilities can occur and in particular the effect of dimension.  It is observed numerically that spherical waves can be unstable depending on the coefficients. A simpler mathematical framework is  to study transversal instability, that means a one dimensional wave propagating in two space dimensions. Then, explicit analytical formulas give explicitely the range of paramaters for instability. 
\end{abstract}

\\[2pt] 
\noindent {\bf Mathematical Classification numbers:} 35C07; 70K50; 76E17; 80A25; 92C17

%%%%%%%%%%%%%%%%%%%%%%%%%%%%%%%%%%%%%%%%%%%%
\section{Introduction}
\label{sec:intro}
%-------------------------------------------
%%%%%%%%%%%%%%%%%%%%%%%%%%%%%%%%%%%%%%%%%%%%
%%%%%%%%%%%%%%%%%%%%%%%%%%%%%%%%%%%%%%%%%%%%

The propagation of unstable interfaces is at the origin of remarkable patterns that can be observed in nature and in experiments. The phenomena has attracted the attention of physicists, geophysicists, chemists and biologists and basic mathematical models can account for this type of unstable dynamical patterns. These models are reaction-diffusion systems and the simplest model is an extension of the Fisher/KPP equation to a two species reaction. It models  reactant concentration and temperature when used for flame propagation \cite{BNS, Billingham}, bacterial population and nutrient concentration when used in biology \cite{mimura, BenJ},  cancer cells and available oxygen/glucosis when used for tumor growth \cite{Benamar_Goriely_2005, sherratt_chaplain, PQV}.

Numerical simulations show that spherical waves can be unstable or stable depending on the model coefficients. But among the many scenarios of instability,  the so-called `transversal instabilities' are the simplest to analyze and explain this surprising effect of dimension which is to de-stabilize a stable one dimensional  traveling wave. The phenomena was observed and related to Diffusion Limited Aggregation, with a first analysis, in \cite{Kessler_Levine,MR1257144}. 

Our goal here is to study such a case of transversal instability and more precisely to understand the modalities of appearance of transversal instabilities for a very simple example given by system. For this, we consider the following two-component reaction-diffusion system~:
\beq
\bepa
\p_t u -  \alpha \Delta u = \f{1}{\alpha} h(u) v,
\\[2mm]
\p_t v - \Delta v = - \f{1}{\alpha} h(u) v .
\eepa
\label{systini}
\eeq
The parameter $\al >0$ is called the Lewis number for flame propagation theory, $\al >1$ is relevant for combustion and $\al <1$ is more relevant for applications to bacterial movement.

Two different cases are proposed both for combustion and biology litteratures depending on propertie of the function $h(\cdot)$
$$
h \in C^1(0,1), \qquad h(0) =0, \qquad h(1)=1,  \qquad h'(u)>0 \quad \text{for }\; 0<u\leq 1, \qquad \text{(KPP type)},
$$
$$ \begin{cases}
h(u) =0 \quad \text{for }\; 0 \leq u \leq \theta<1, \qquad h(\theta+)=h^+\geq 0, 
\\[2mm]
h(1)=1, \quad h'(u) \geq 0 \quad \text{for }\; \theta <u\leq 1,
\end{cases}
\qquad \text{(Ignition temperature type).}
$$
Our interest lies on two dimensional stability of one dimensional traveling waves for this system. A proof of existence for one dimensional traveling wave solutions can be found in \cite{marion_85} when $h(\cdot)$ is of KPP type and when  $\al \geq 1$ for ignition temperature type. Also, in \cite{BNS}, the authors  prove existence of traveling waves when $h(\cdot)$ is of ignition temperature type and no restriction on $\al$. More recent results for KPP type, in a cylinder and covering all Lewis numbers, can be found in \cite{hamel_ryzhik}.

In this paper, we consider a simple example for which we can handle analytical 
computation. It corresponds to ignition temperature type and the function $h$ 
is given by 
\beq\label{eq:h}
h(u) = \left\{
\begin{array}{l} 0, \qquad \mbox{for } u\leq \theta, \\
1,\qquad \mbox{for } u>\theta. \end{array}
\right.
\eeq
We first report in Section~\ref{sec:simul}, based on numerical simulations, two dimensional spherical waves which can be unstable for certain coefficients. Then,  we build analytically the one dimensional traveling waves in Section~\ref{sec:tw}. The analytical formulas are fundamental to handle the spectral problem arising to study linearized stability of the transversal waves in Section~\ref{sec:stability}.

%-------------------------------------------
\section{Numerical observations}
\label{sec:simul}
%-------------------------------------------

In two dimensions, numerical simulations of system \eqref{systini}--\eqref{eq:h} 
exhibit various behaviours depending on the values of $\alpha$ and $\theta$ 
in $(0,1)$. We present them here as a motivation for our theoretical study. 

These simulations are obtained using the finite element method implemented within the software FreeFem++  \cite{FreeFem, hecht}. 
The computational domain is a disc with radius $4$ and we denote by $\Gamma$ its
boundary. At the boundary, Neumann boundary conditions are implemented 
for both $u$ and $v$~:
 $$
 \p_\nu u|_{\Gamma}=0,\qquad \p_\nu v|_{\Gamma}=0,
 $$ 
where $\nu$ is the outward unit normal. 
We use a semi-implicit time discretization.
Then the resulting system is discretized thanks to P1 finite element method.

\begin{figure}[ht]
\centerline{\includegraphics[width=7cm]{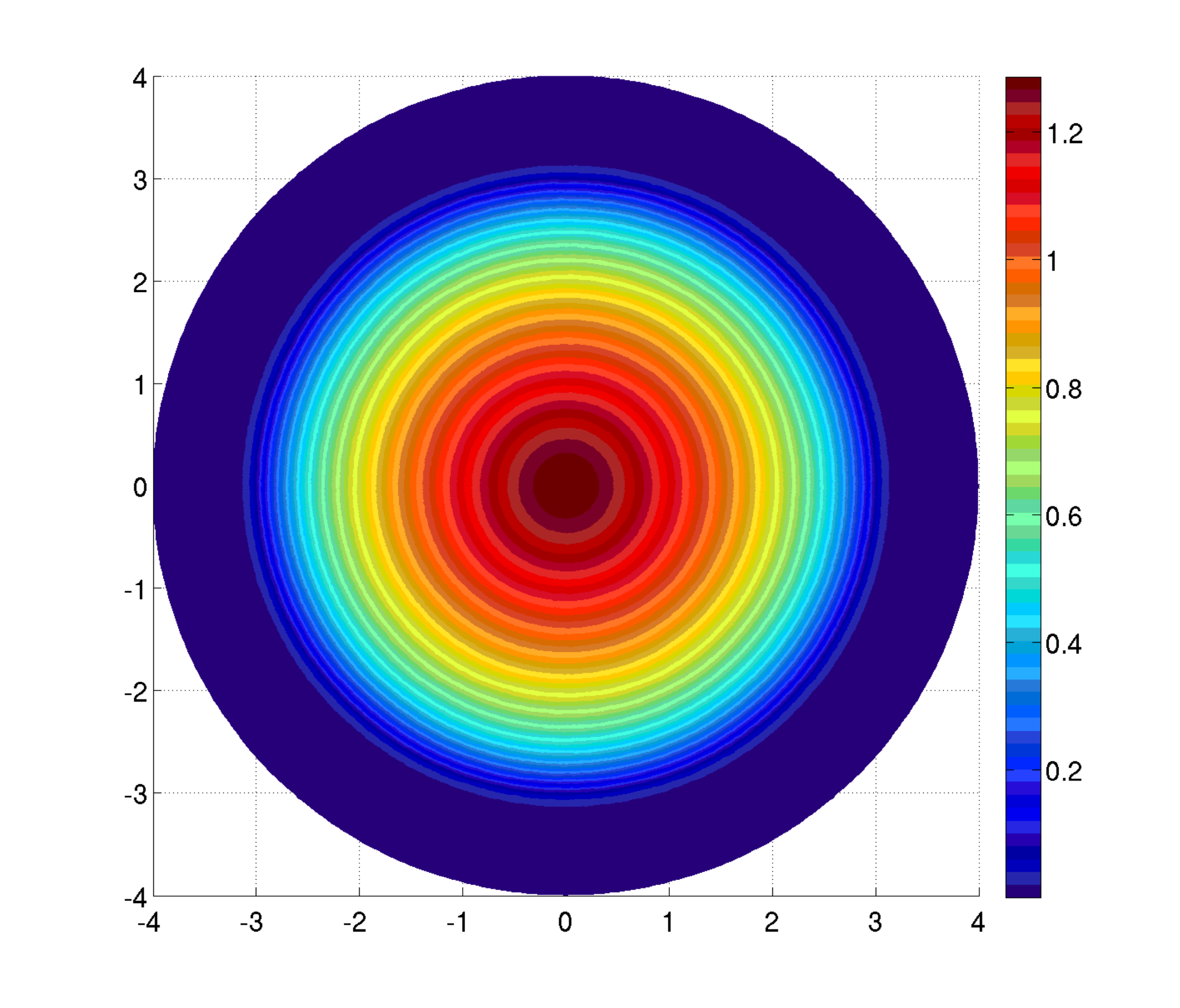}
\hspace{10pt}\includegraphics[width=7cm]{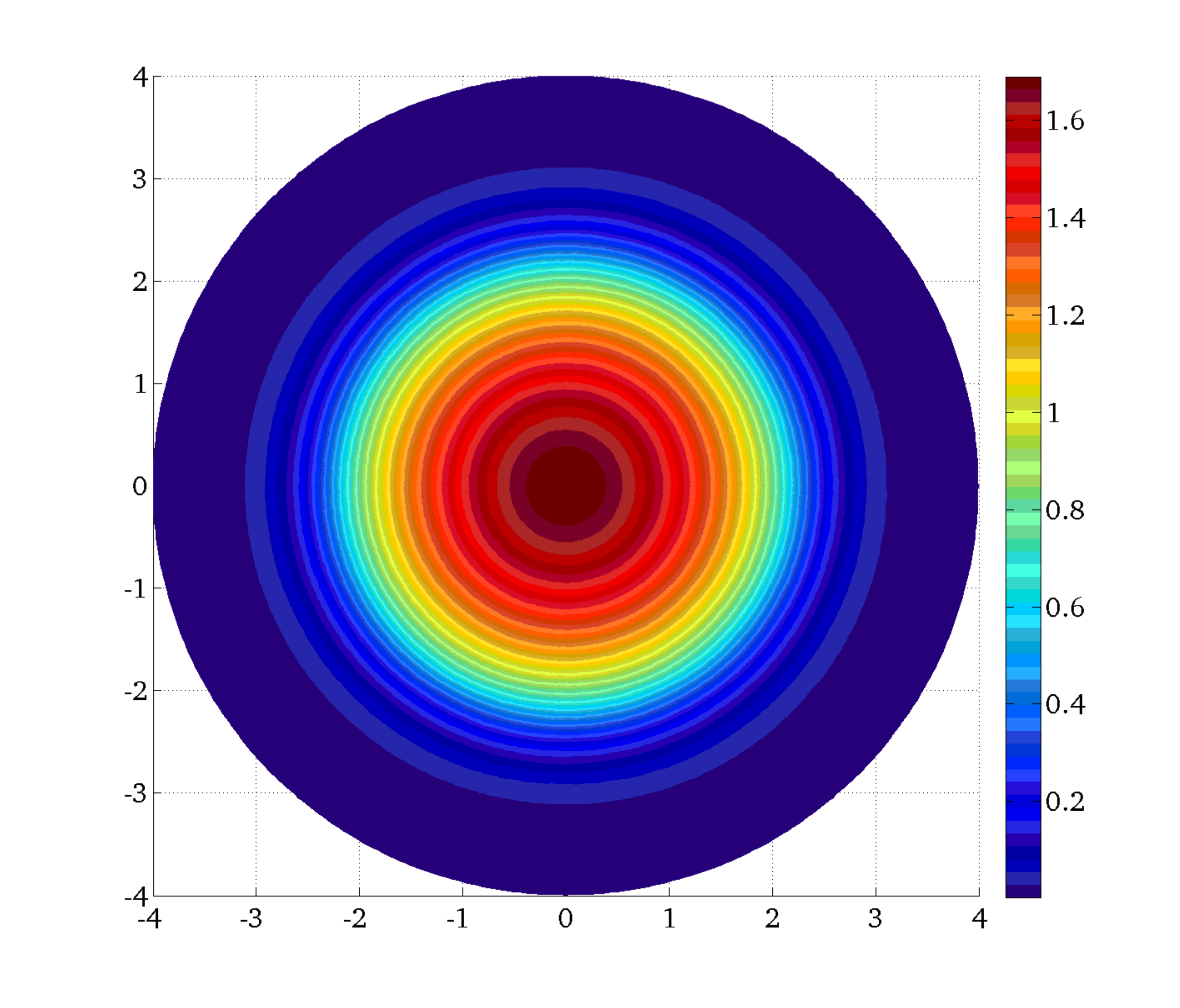}}
\caption{\small\sl\label{simul1} Numerical simulations for component $u$ in \eqref{systini}. (1) Left: Plot of the computed $u$ at time $T=1$ for $\alpha = 0.25$ and $\theta=0.1$. (2) Right: Plot of the computed $u$ at time $T=3$ for $\alpha = 0.25$
and $\theta = 0.5$.
In both cases, we observe the propagation of a stable spherical wave which invades the computational domain.}
\end{figure}

\begin{figure}[ht]
\centerline{\includegraphics[width=7cm]{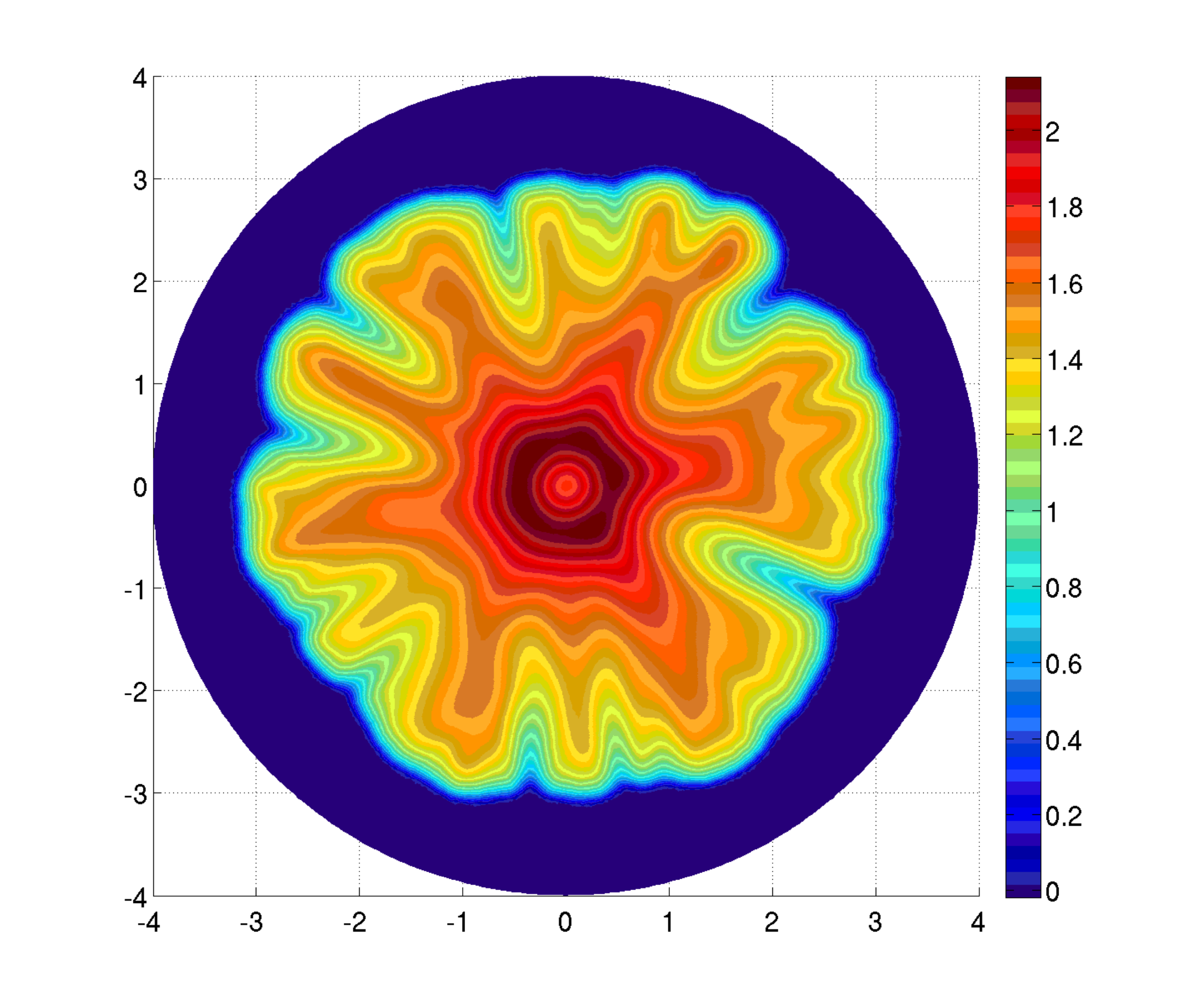}
\hspace{10pt}\includegraphics[width=7cm]{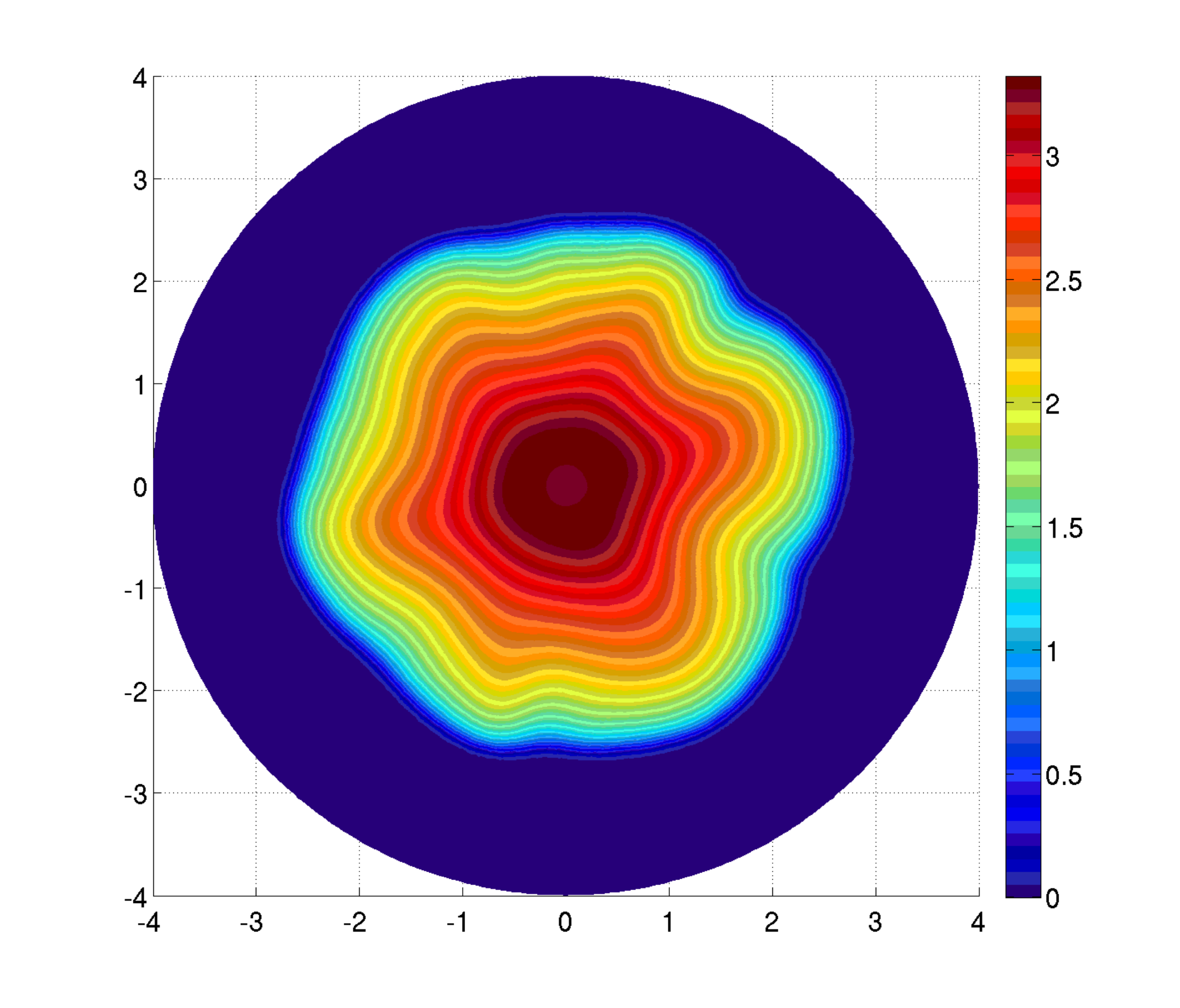}}
\caption{\small\sl\label{simul2} Numerical simulations for species $u$ in \eqref{systini}. (1) Left: Plot of the computed $u$ at time $T=2$ for $\alpha = 0.01$ and $\theta=0.1$. (2) Right: Plot of the computed $u$ at time $T=6$ for $\alpha = 0.25$
and $\theta = 0.5$.
When $\alpha$ is small, numerical instabilities appear which create a complex fingered pattern.}
\end{figure}

The initial given data is~:
$$
u^0 = {\bf 1}_{\{\sqrt{x^2+y^2}\leq 0.4\}}; \qquad v^0 = 1-u^0.
$$
The time step is $dt=0.0025$ and the number of nodes is $21879$.
The numerical results are depicted in Figure \ref{simul1}
and \ref{simul2} where the computed approximation of $u$
is plotted after several time iterations for different values
of the parameters $\alpha$ and $\theta$.
Depending on the values of $\alpha$ and $\theta$, 
we observe different patterns in the numerical simulations.
Figure \ref{simul1} displays the numerical simulations for $\alpha = 0.25$
and for $\theta=0.1$ (Left) and $\theta=0.5$ (Right). 
In both cases, we observe a wave that invades the computational domain 
and the numerical result do not show instabilities.
Comparing this two results, we deduce that the invasion process depends on 
$\theta$.
Figure \ref{simul2} displays the numerical results
obtained for small $\alpha$~: we choose $\alpha = 0.01$. 
In this case, we observe numerical instabilities that create a complex pattern.
Instabilities are much more visible when $\theta = 0.1$ (Figure \ref{simul2},
Left) than when $\theta = 0.5$ (Figure \ref{simul2}, Right).

%-------------------------------------------
\section{One dimensional traveling waves}
\label{sec:tw}
%-------------------------------------------

One dimensional  traveling waves are solutions of the form
$$
u(t,x)=u_0(x-\sigma t), \quad v(t,x)=v_0(x-\sigma t),
$$
where $\sigma>0$ is a constant representing the traveling wave velocity.  They are a convenient way to understand the propagation phenomena presented in Section~\ref{sec:simul}.

For system \eqref{systini} traveling waves  are determined from the system~:
\beq
\bepa
- \sigma {u_0}_x -  \alpha {u_0}_{xx}= \f{1}{\alpha} h(u_0) v_0,
\\[2mm]
- \sigma {v_0}_x - {v_0}_{xx}= - \f{1}{\alpha} h(u_0) v_0,
\\[2mm]
(u_0,v_0)(-\infty)= (1,0), \qquad  (u_0,v_0)(+ \infty)= (0,1).
\eepa
\label{eq:uv}
\eeq
To avoid ambiguity due to the translation invariance of  the problem we set
\beq
0< u_0(0)= \theta <1.
\label{eq:norm}
\eeq
We say that a traveling wave solution to (\ref{eq:uv}) is {\em monotonic}, if each component is monotonic, and then we can  normalize the signs with $u_0'<0$ and $v_0' >0$. 

%-------------------------------------------
\begin{proposition}\label{prop:tw}
There exists {{a unique}} monotonic traveling wave for system \eqref{systini}--\eqref{eq:h},
i.e. {{a unique}} $\sigma>0$ and a pair $(u_0,v_0)\in \mathcal C^{1,\nu}(\R)$, $\nu\in (0,1)$, solving \eqref{eq:uv}--\eqref{eq:norm} 
with $u_0$ nonincreasing,  $v_0$ nondecreasing.

More precisely this traveling wave solution moves with the speed
\beq\label{eq:sigma}
\sigma = (1-\theta)\sqrt{\f{\alpha}{\theta^2+\alpha\theta(1-\theta)}}
\eeq
and  is given explicitly by
\beq\label{eq:u0}
u_0(x)=
\left\{\begin{array}{ll}
1-(1-\theta)e^{\theta x/\beta},  \qquad &\mbox{ for } x<0, \\[2mm]
\theta e^{-(1-\theta)x/\beta},   \qquad& \mbox{ for } x>0, 
\end{array} \right.
\eeq
\beq\label{eq:v0}
v_0(x)=
\left\{\begin{array}{ll}
\f{\alpha(1-\theta)}{\theta+\alpha (1-\theta)} e^{\theta x/\beta},
\qquad & \mbox{ for }x<0, \\[2mm]
1-\f{\theta}{\theta+\alpha (1-\theta)} e^{-\alpha(1-\theta)x/\beta},
\qquad & \mbox{ for }x>0,
\end{array} \right.
\eeq
where
\beq\label{nota:beta}
\beta = \sqrt{\alpha\theta} \,\sqrt{\theta+\alpha(1-\theta)}.
\eeq
\end{proposition}
%----------------------------------------------

\begin{figure}[ht]
\centerline{\includegraphics[width=8cm]{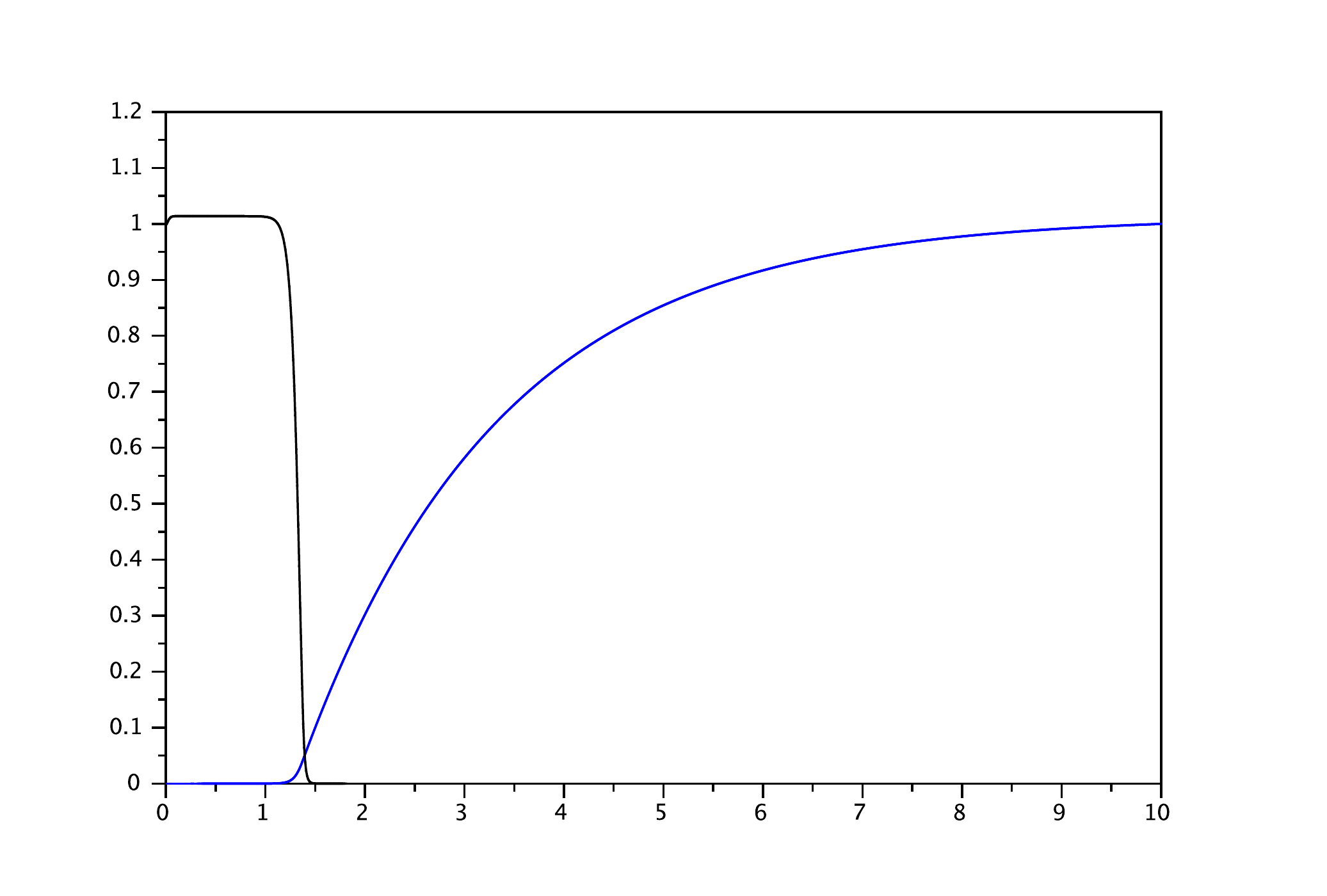}}
\vspace{-7mm}
\caption{\small \sl  The traveling wave solution with $u$ decreasing and $v$ increasing for $\al=0.005$ and $\theta =0.1$.}
\label{Fig:tw}
\end{figure}

This solution is depicted in Figure~\ref{Fig:tw}. 
\begin{remark}
We notice that when $\theta$ goes to $0$, we have that $\sg\sim \f{\sqrt{\alpha}}{\theta}$.
Then the wave goes faster when $\theta$ is smaller. 
This remark confirms our observation in Figure \ref{simul1}, where we can
notice than the invasion process of species $u$ is faster when $\theta$ is smaller.
\end{remark}

%-----------------------------------------------------
%---------------------------------------------------------------------------------------

\begin{proof}
We recall that we look for a nonincreasing function $u$ and we have denoted
$u_0(0)=\theta$. Therefore, from the definition of the nonlinearity  $h(\cdot) $ in \eqref{eq:h},
we have $h\big(u_0(x)\big)=0$ for $x>0$. Then, the system \eqref{eq:uv} is reduced to 
$$
- \sigma {u_0}_x - \alpha {u_0}_{xx}= 0, \qquad
- \sigma {v_0}_x - {v_0}_{xx}= 0, \qquad x >0.
$$
With the boundary condition at $+\infty$ in \eqref{eq:uv}~: $(u_0,v_0)(+\infty)=(0,1)$, we deduce
\beq\label{eq:uvpos}
u_0(x)=\theta e^{-\sg x /\alpha}, \qquad 
v_0(x) = 1-be^{-\sg x}, \qquad \mbox{ for } x>0,
\eeq
where $b$ is a constant to be fixed later.

For $x<0$, we have $h\big(u_0(x)\big)=1$, therefore system \eqref{eq:uv} is reduced to
$$
- \sigma {u_0}_x - \alpha {u_0}_{xx}= \f{v_0}{\alpha}, \qquad
\sigma {v_0}_x + {v_0}_{xx}- \f{1}{\alpha} v_0= 0.
$$
Solving the second equation, and because the solution is continuous at $x=0$ by ellipitic regularity, leads to
\beq\label{eq:vneg}
v_0(x)=(1-b)e^{\lambda x}, \qquad \mbox{ with }\quad 
\lambda = \f 12 \big(-\sigma+\sqrt{\sg^2+4/\alpha}\big),
\qquad \mbox{ for } x<0,
\eeq
where we have used the boundary conditions~: $(u_0,v_0)(-\infty)=(1,0)$.
Then, we obtain
\beq\label{eq:uneg}
u_0(x) = 1-(1-\theta)e^{\lambda x},\qquad \mbox{ for } x<0,
\eeq
which is a solution of the equation for $u_0$ provided 
$$
\big(\sigma \lambda +\alpha \lambda^2\big)(1-\theta)= \f{1-b}{\alpha}.
$$
This latter equality allows to determine the value of $b$~:
$$
b=1-\alpha\big(1+\sg \lambda(1-\alpha)\big) (1-\theta).
$$
Finally, the continuity of the derivative $u'_0(0^+)=u'_0(0^-)$ implies
$\dis \f{\sg\theta}\alpha = \lambda (1-\theta).$
Using this relation and  the expression of $\lambda$  \eqref{eq:vneg} we obtain
$$
\sigma = (1-\theta)\sqrt{\f{\alpha}{\theta^2+\alpha\theta(1-\theta)}}.
$$
We deduce
\beq\label{eq:lambdab}
\lambda = \frac{\sqrt{\theta}}{\sqrt{\alpha}\sqrt{\theta+\alpha(1-\theta)}}, 
\qquad \sigma\lambda = \frac{1-\theta}{\theta+\alpha(1-\theta)}, 
\qquad b = \f{\theta}{\theta+\alpha(1-\theta)}.
\eeq

The conclusions stated in Proposition~\ref{prop:tw} follow directly from the construction and  formulas  \eqref{eq:uvpos}--\eqref{eq:vneg}--\eqref{eq:uneg}.
\end{proof}
\qed

%%%%%%%%%%%%%%%%%%%%%%%%%%%%%%%%%%%%%%%%%%%%
\section{Stability of planar traveling waves}
\label{sec:stability}
%-------------------------------------------
%%%%%%%%%%%%%%%%%%%%%%%%%%%%%%%%%%%%%%%%%%%%
%%%%%%%%%%%%%%%%%%%%%%%%%%%%%%%%%%%%%%%%%%%%

As suggested by the numerical results in Section \ref{sec:simul}, based on spherical waves, we expect that 
transversal instability can occur in two dimensions.

We propose here to study the linear transversal stability. To do so, and in the spirit of \cite{Kessler_Levine, BenAmar_C_F} for instance, with $\e\ll 1$,  we set
$$
\bepa
u(t, x,y) = u_0(x-\sigma t) + \e e^{\lambda t} \cos(\omega y) u_1(x-\sigma t),  \\[2mm]
v(t, x,y) = v_0(x-\sigma t) + \e e^{\lambda t} \cos(\omega y) v_1(x-\sigma t).
\eepa
$$

Substituting this expansion into \eqref{systini} and keeping only the term 
of order 1 in $\e$, we get the linearized system (in the traveling wave frame)
\beq
\bepa
\lambda u_1 - \sigma u_1' -\alpha u_1'' + \alpha \omega^2 u_1 = \f{1}{\alpha}
\big(h'(u_0)v_0 u_1 + h(u_0) v_1\big),   \\[2mm]
\lambda v_1 - \sigma v'_1 - v_1'' + \omega^2 v_1 = -\f{1}{\alpha}
\big(h'(u_0)v_0 u_1 + h(u_0) v_1\big).
\eepa
\label{eq:u1v1}
\eeq
We notice that for $\omega =0$, the system has the solution $\lambda=0$, 
$u_1=u'_0$ and $v_1=v'_0$. Notice that $\omega=0$ also  represents the case of dimension one (no transversal effect)s.

\begin{definition}\label{def:unstable}
In two dimensions, we say that the one dimensional traveling wave $(u_0,v_0)$ 
for system \eqref{systini} in Proposition~\ref{prop:tw} is 
{\bf transversally linearly unstable}
if there exists $\omega >0$ and $\lambda$ with $\mathrm{Re}\, \lambda>0$ such that 
system \eqref{eq:u1v1} admits a non-trivial  solution in $\mathcal C^\nu(\R)\cap L^2(\R)$.
\end{definition}

\begin{proposition}\label{prop:instab}
Let $\theta\in (0,1)$. Let us consider the function $h$ given in \eqref{eq:h}.
Then the following hold:
\begin{enumerate}
\item  For $\alpha$ small enough,  the traveling waves in Proposition \ref{prop:tw} are linearly stable in one dimension for all $\theta \in (0, 1)$.

\item For each $\theta\in (0,1)$ and each small $\alpha$ there exists $\omega(\alpha,\theta)$ such that $\lambda(\omega(\alpha, \theta))<0$ i.e. the traveling wave is transversally linearly unstable for these values of the parameters in two dimensions. Moreover $\omega(\alpha,\theta)=\mathcal O(\frac{1}{\sqrt{\alpha}})$. 
\end{enumerate}
\end{proposition}

The paper \cite{sattinger_system} suggests that, in appropriate weighted spaces, one dimensional traveling wave are nonlinearly stable in the range of parameters when they are linearly stable. 
\\

\begin{proof}
Since linear stability in one dimension reduces to studying (\ref{eq:u1v1}) for $\omega=0$, from now on we assume more generally that $\omega\geq 0$. We will  look for  $\lambda\in \R$, $\lambda=\lambda(\omega)>0$
such that system \eqref{eq:u1v1} admits a non-trivial solution. Note that $\lambda$ depends on $\alpha$ and $\theta$ as well, but we will not make this dependence explicit unless necessary. 

For $x>0$, system \eqref{eq:u1v1} reduces to 
$$
\bepa
(\alpha \omega^2+\lambda) u_1 - \f{(1-\theta)\alpha}{\beta} u'_1-\alpha u''_1 =0, \\[2mm]
(\lambda+\omega^2) v_1 - \f{(1-\theta)\alpha}{\beta} v'_1-v''_1=0.
\eepa
$$
We can solve this linear problem and obtain
\beq\label{u1v1pos}
\bepa
u_1(x)=A e^{r_- x}, \qquad r_- = -\f{(1-\theta)}{2\beta}-\f{1}{2\beta}\sqrt{(1-\theta)^2+4\beta^2(\omega^2+\f{\lambda}{\alpha})}, \\[2mm]
v_1(x)= B e^{s_- x}, \qquad s_- = -\f{(1-\theta)\alpha}{2\beta}-\f{1}{2\beta}\sqrt{(1-\theta)^2\alpha^2+4\beta^2(\omega^2+\lambda)}.
\eepa
\eeq
Here,   $A$ and $B$ are constants to be determined and $r_\pm$ are the roots of the polynomial
\beq
 (\alpha \omega^2+\lambda) - \f{(1-\theta)\alpha}{\beta} r -  \al r^2 =0.
\label{magic_pol}
\eeq

For $x<0$, system \eqref{eq:u1v1} reduces to
$$
\bepa
(\alpha \omega^2+\lambda) u_1 - \f{(1-\theta)\alpha}{\beta} u'_1-\alpha u''_1 =\f{1}{\alpha} v_1, \\[2mm]
\Big(\lambda+\omega^2+\f{1}{\alpha}\Big) v_1 - \f{(1-\theta)\alpha}{\beta} v'_1-v''_1=0,
\eepa
$$
where we have used the expression of $\sigma$ in \eqref{eq:sigma} recalling
that $\beta$ is given in \eqref{nota:beta}.
Then we get
\beq\label{eq:mu-}
v_1(x)= B e^{\mu_+ x}, \qquad \mu_+ = -\f{(1-\theta)\alpha}{2\beta}  + \f{1}{2\beta}\sqrt{(2\theta+\alpha(1-\theta))^2+4\beta^2(\omega^2+\lambda)}.
\eeq
Substituting  this expression in the equation for $u_1$, we get
\beq\label{eq:r+}
u_1(x)=C e^{r_+ x}+\gamma B e^{\mu_+ x}, \qquad r_+ = -\f{(1-\theta)}{2\beta}+\f{1}{2\beta}\sqrt{(1-\theta)^2+4\beta^2(\omega^2+\f{\lambda}{\alpha})}, \\[2mm]
\eeq
where $C$ is a constant and we have set
$$
\Big(\alpha \omega^2+\lambda - \f{(1-\theta)\alpha}{\beta} \mu_+-\alpha \mu_+^2 \Big)
\gamma = \f{1}{\alpha}.
$$
The value of the parameter $\gamma$ follows from the definition of $r_\pm$ as the roots of \eqref{magic_pol}
\beq\label{eq:gamma}
\gamma = -\f{1}{\alpha^2(\mu_+-r_+)(\mu_+-r_-)}.
\eeq
Moreover, by continuity at $x=0$, we need $C+B\gamma = A$.
By definition of the function $h$ and with \eqref{eq:u0}--\eqref{eq:v0}, we obtain
\beq\label{eq:delta_h}
h'(u_0)v_0 = \f{\beta \alpha}{\theta(\theta+\alpha(1-\theta))}\delta_{x=0}
=\f{\alpha^2}{\beta} \delta_{x=0}.
\eeq
As a consequence, the jump relation for $(u'_1,v_1')$ at $x=0$, which  can be  deduced from equation \eqref{eq:u1v1}, leads to
$$
\bepa
-\alpha(u'_1(0^+)-u'_1(0^-)) = \f{\alpha}{\beta} u_1(0), \\[2mm]
-v'_1(0^+)+v'_1(0^-)=-\f{\alpha}{\beta} u_1(0).
\eepa
$$
Writing these equalities in terms of the free parameters, we arrive to the following set of  relations
$$
\bepa
A=C+\gamma B,
 \\[2mm]
A\big(r_-+\f{1}{\beta}\big) = C r_+ + \gamma \mu_+ B, 
\\[2mm]
\f{\alpha}{\beta} A = B s_- - B \mu_+. 
\eepa
$$
Replacing $B$ and $C$ in the second equation, we get
$$
A\big(r_-+\f{1}{\beta}\big) = A r_+ + A \gamma (\mu_+-r_+) \f{\alpha}{\beta(s_--\mu_+)}.
$$
We conclude that there exists a non trivial solution to \eqref{eq:u1v1}
provided the following identity holds~:
\beq\label{relation}
1 = \beta(r_+-r_-) + \f{1}{\alpha(\mu_+-r_-)(\mu_+-s_-)},
\eeq
where we used the expression \eqref{eq:gamma}.
We verify straightforwardly that for $\omega=0$ and $\lambda =0$ 
relation \eqref{relation} is always satisfied.

It remains to compute the value of $\lambda$ using \eqref{relation}. In this algebraic  equation  $\lambda$ is given implicitly as a function of the parameters $\alpha, \theta$ and $\omega$, and in the analysis of this expression we rely on taking the limit $\alpha\to 0$ and also on Maple based simulations (in this sense our proof is to some small extent computer assisted).
 
First, we consider the case $\omega=0$ which provides the stability in one dimension. Then we consider the limit $\alpha\to 0$ choosing the scale: $\omega = \frac{\omega_0}{\sqrt{\alpha}}$, with $\omega_0>0$ fixed.

\noindent
{\em Case 1: $\alpha\to 0$, $\omega=0$.} This case covers Assertion 1 of the propostion. 
We set 
$$
\zeta=\sqrt{(1-\theta)^2+4\beta^2\lambda/\alpha}, \qquad \eta = \sqrt{(2\theta+\al(1-\theta))^2+4\beta^2\lambda}.
$$
Using \eqref{u1v1pos}, \eqref{eq:mu-}, \eqref{eq:r+}, the identity 
\eqref{relation} in the case $\omega=0$ reduces to \begin{equation}
1 = \zeta+ \f{4\beta^2}{\al \big((1-\theta)(1-\al)+\zeta+\eta \big)(\eta+\al \zeta)}.
\label{cold fact}
\end{equation}
%\textcolor{blue}{It is not difficult to check numerically that 
%the right hand side is increasing with respect to $\lambda$ when 
%$\lambda>0$. But I do not succeed in proving it. If everything goes right for $\lb=0$ is should be larger than 1}.
Taking $\alpha=0$ above, we find that 
\[
\lambda=-\frac{1}{2} \, \frac{1-\theta}{\theta}.
\]
By continuity, this means in particular that for all sufficiently small $\alpha$, the traveling wave solution is linearly stable. Hence, it is tempting to speculate that in fact linear stability is true for any $\alpha>0$, see however Remark~\ref{sixto rodriguez}.

%\medskip
%\noindent {\underline{Case 2: $\alpha \to 0$ and $\omega$ bounded.}}
%We use formula (\ref{relation}). In the limit $\alpha\to 0$ we get
%\[
%\beta(r_+-r_-)=\beta\left[\frac{1}{\beta}\sqrt{(1-\theta)^2+4\beta^2\left(\omega^2+\frac{\lambda}{\alpha}\right)}\right]\longrightarrow \sqrt{(1-\theta)^2+ 4\theta^2\lambda},
%\]
%where we have used (\ref{nota:beta}).  Next, using (\ref{u1v1pos}) and (\ref{eq:mu-}) we get
%\[
%\begin{aligned}
%\sqrt{\alpha}(\mu_+-r_-)&\longrightarrow 1-\frac{1-\theta}{2\theta}-\frac{1}{2\theta}\sqrt{(1-\theta)^2+4\theta^2\lambda}, \\
%\sqrt{\alpha}(\mu_+-s_-)&\longrightarrow 1.
%\end{aligned}
%\]
%Denote $\zeta=\sqrt{(1-\theta)^2+4\theta^2\lambda}$. Then we need to solve for $\zeta\geq 0$ the equation
%\[
%1-\zeta=\frac{2\theta}{1+\theta+\zeta}\Longrightarrow \zeta=\frac{-\theta +\sqrt{\theta^2+4(1-\theta)}}{2}
%\]
%from which we find
%\[
%\lambda=\lambda(\theta) = \frac{\zeta^2}{4\theta^2}-\frac{1}{4}\left(\frac{1}{\theta}-1\right)^2=0.
%\]
%This means that  oscillations with relatively large or relatively small (with respect to ${\sqrt{\alpha}}$, see  case 2 below)  period $\frac{2\pi}{\omega}$ will {\it not} give rise to instability.  
%
\medskip

\begin{figure}[ht]
\centerline{\includegraphics[width=6cm]{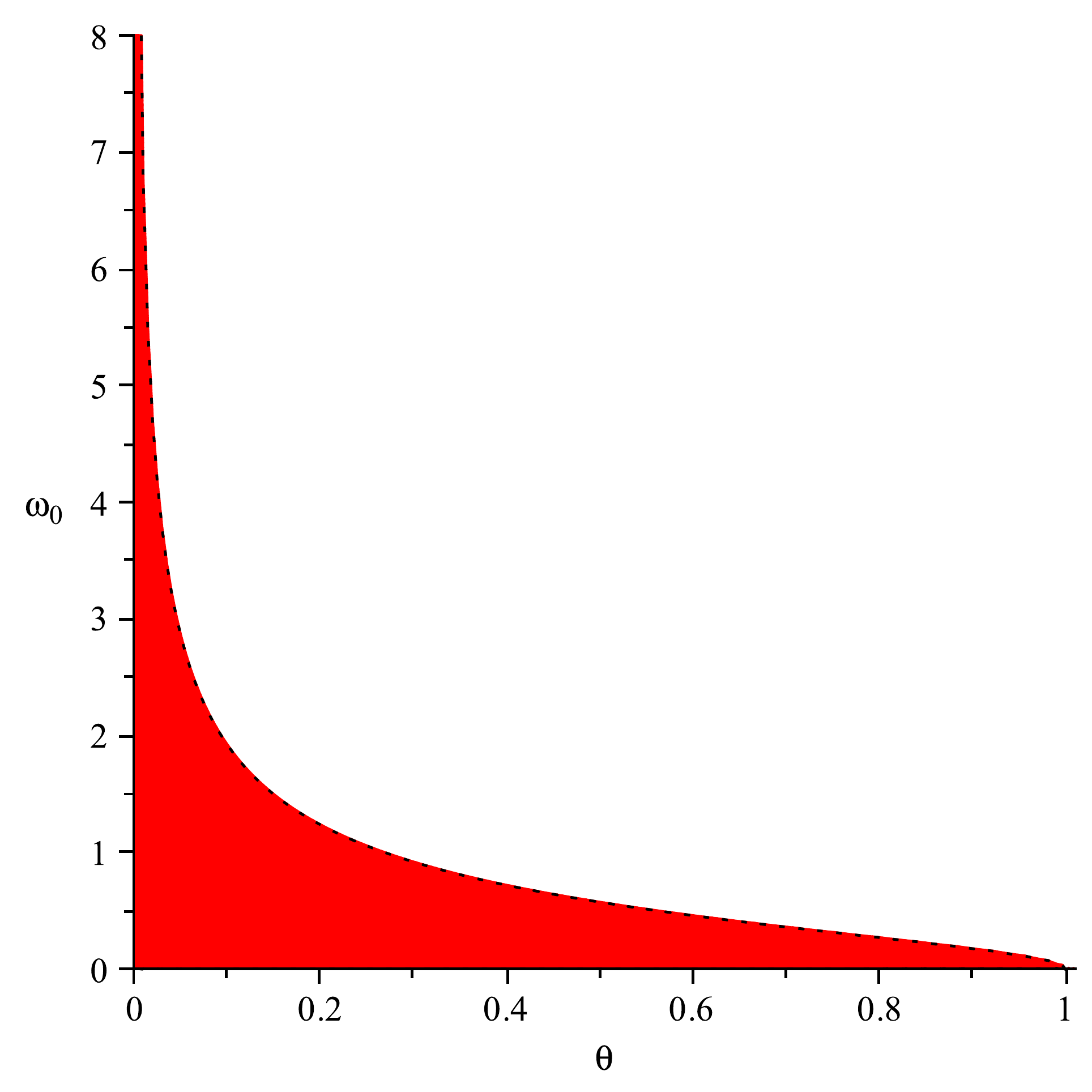}}
\vspace{-7mm}
\caption{\small \sl  Plot of the region of instability  $\lambda=\lambda(\theta, \omega_0)>0$, marked red.}
\label{Fig1}
\end{figure}

\noindent
{\em Case 2: $\alpha\to 0$ and $\omega=\frac{\omega_0}{\sqrt{\alpha}}$.}
Calculations are quite similar in this case. Denoting now
\[
\zeta=\sqrt{(1-\theta)^2+4\theta^2(\lambda+\omega^2_0)},
\]
we need to solve
$$
1-\zeta = \frac{2\theta}{(\sqrt{\omega_0^2+1}+\omega_0)(1-\theta+2\theta\sqrt{\omega_0^2}+\zeta)}.
$$
Therefore,
$$
\zeta=\frac{1}{2}\left(\theta-2\theta\sqrt{\omega_0^2+1}
+\sqrt{\Big(\theta-2\theta\sqrt{\omega_0^2+1}\Big)^2+4-4\theta+8\theta\omega_0}\right)
$$
and then
$$
\lambda = -\omega_0^2 + \frac{1}{16\theta^2} \Big(
\theta-2\theta\sqrt{\omega_0^2+1}+\sqrt{(\theta-2\theta\sqrt{\omega_0^2+1})^2
+4-4\theta+8\theta\omega_0}\Big)^2-\Big(\frac{1-\theta}{2\theta}\Big)^2.
$$
In this case, depending on the value of $\omega_0$ we may have $\lambda>0$ or $\lambda<0$. Figure~\ref{Fig1} illustrates the situation. 

\end{proof}

%-----------------------------------------------------
%---------------------------------------------------------------------------------------

\qed

\section{Concluding remarks}
\label{judas priest}

\begin{remark}
It has been noted with the numerical simulations of Section\ref{sec:simul},
that for small values of $\alpha$, instabilities are more visible when 
$\theta$ is small (see Figure \ref{simul2}). 
The plot of the region of instability in Figure \ref{Fig1} confirms this observation.
In fact, we notice on this latter Figure that small values of $\theta$
allows large values of $\omega_0$ which can be seen as a frequence of oscillations in the transversal direction.
\end{remark}

\begin{remark}\label{sixto rodriguez}
We use formula (\ref{cold fact}) and denote
\[
F(\alpha, \theta, \lambda)=\zeta+ \f{4\beta^2}{\al \big((1-\theta)(1-\al)+\zeta+\eta \big)(\eta+\al \zeta)}-1.
\]
 Solutions of $F(\alpha, \theta, \lambda)=0$ determine the eigenvalues  $\lambda=\lambda(\alpha,\theta)$. It is easy to check that 
\[
\lim_{\alpha\to \infty} F(\alpha, \theta, \lambda)=\infty, \quad \forall\lambda>0, \; \theta\in(0,1),
\]
however this limit is not uniform. Indeed there exist a  $\theta^*$ such that for any $\theta \in(\theta^*,1)$ we can find a value $\alpha (\theta) >0$ such that $F(\alpha, \theta, \lambda)=0$ for some $\lambda>0$. We illustrate this in  Figure~\ref{wicked lady}.
%-------------------------------------------------------------
\begin{figure}[ht]
\centerline{\includegraphics[width=7cm]{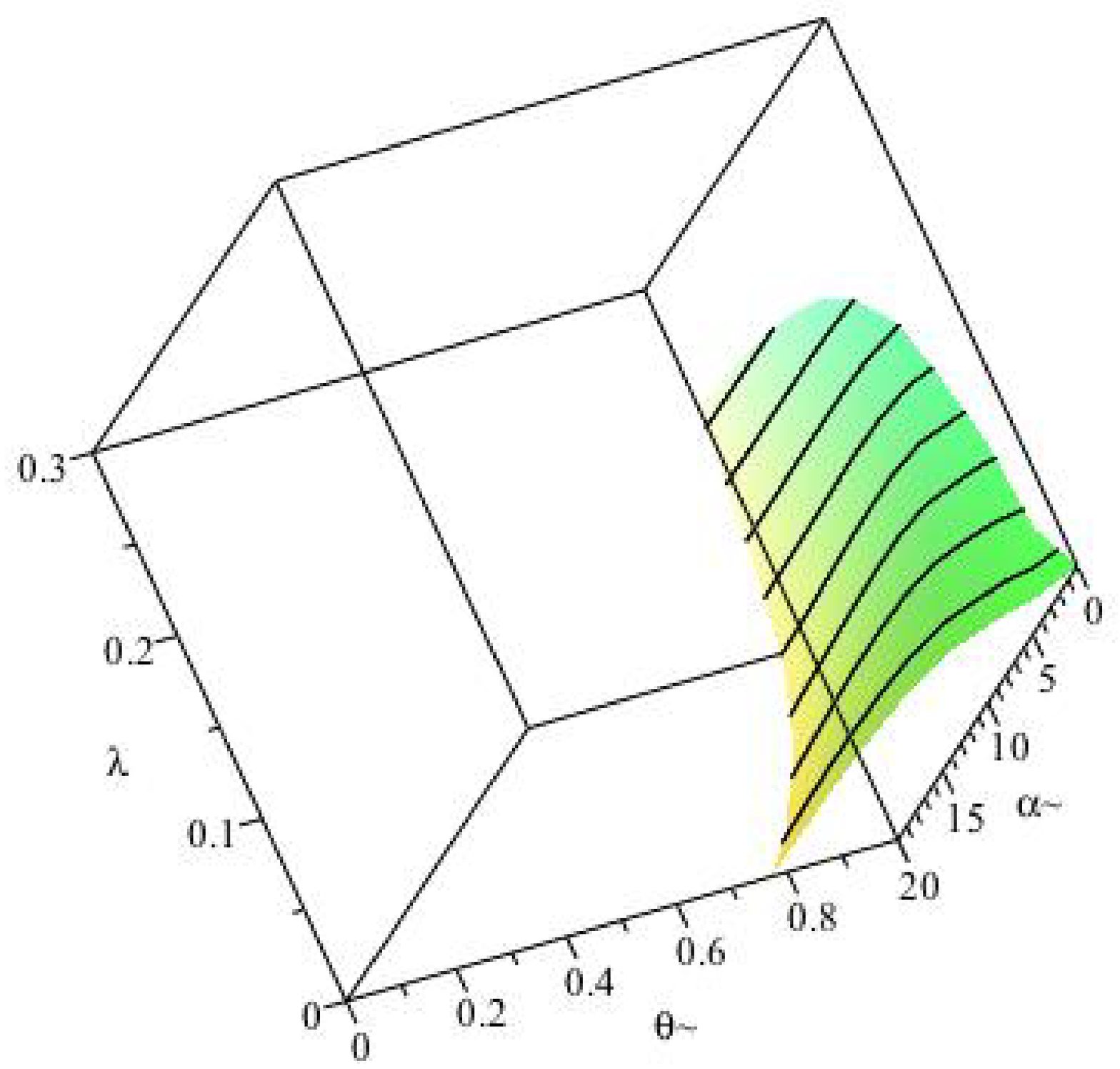}
\hspace{10pt}\includegraphics[width=7cm]{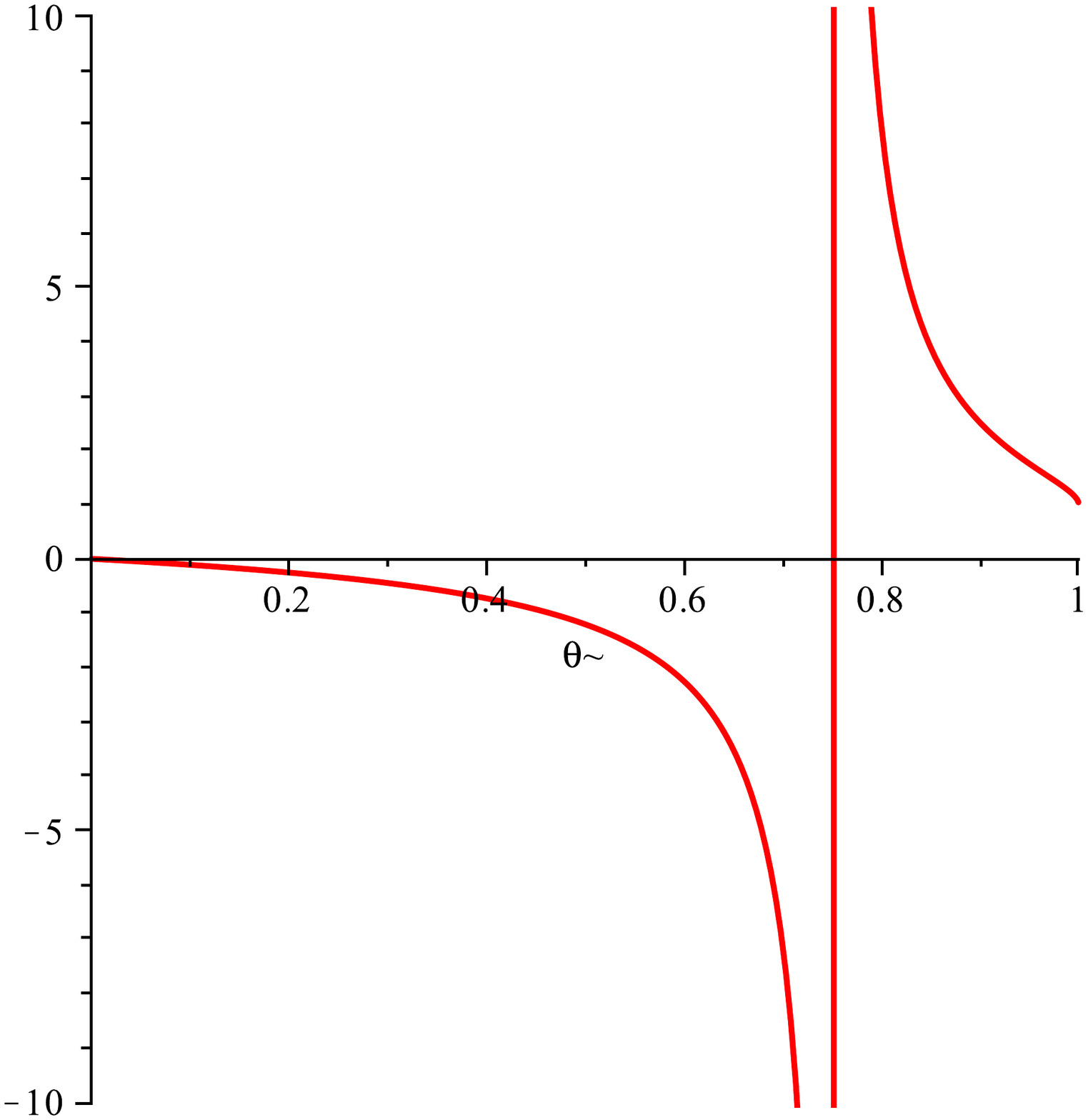}}
\vspace{-12mm}
\caption{\small \sl  (Left) Plot of the level set $F(\alpha, \theta, \lambda)=0$. (Right) Plot of the implicitly defined  curve $\theta\mapsto \alpha(\theta)$ where $F(\alpha(\theta), \theta, 0)=0$. It suggests that the instability appears for some $\theta^*\in (0.7, 0.8)$ and the corresponding values of $\al$ are larger than $1$. This is confirmed by a more refined analysis of the picture on the right.}
%\caption{\small\sl\label{Fig1}  (1) Left: Plot of the region $\lambda=\lambda(\theta, \omega_0)<0$ (green) (2) Right: Plot of  the curve $\lambda=\lambda(\theta; \omega_0)$ for various values of $\theta$ with $\omega_0$ varying.}
\label{wicked lady}
\end{figure}

The fact that the traveling wave in one dimension is unstable for large values of the Lewis number is somewhat of a surprise and raises a more general question of stability or instability for the problems with KPP type or ignition type nonlinearities. Note that, in both cases, we are dealing with a prey-predator system and in particular the linear problem is non-cooperative. This means that   methods based on maximum principle do not work and the known results (see for instance \cite{MR1297766})) do not apply. The special feature of our problem is the monotonicity of the traveling fronts. With this property, one may expect that they should be stable, as it happens for scalar problems and is seen easily there  from the Krein-Rutman theorem. For systems of equations, there is no general theory that one could apply but monotone waves are stable in some cases (see for instance \cite{MR3062741}). Our example shows that the question is in fact more subtle.
\end{remark}

\begin{remark}\label{buzzcocks}
Writing (\ref{relation}) in the form 
\[
0=-1+ \beta(r_+-r_-) + \f{1}{\alpha(\mu_+-r_-)(\mu_+-s_-)},
\]
we obtain a dispersion relation $\lambda=\lambda(\omega)$ for any fixed $\alpha$ and $\theta$. Pictures in Figure~\ref{ever fall in love}, confirm the intuitively obvious fact that there should always be the most unstable frequency $\omega^*$, that is a maximum value of $\lambda$ and this is relevant of Turing instability.
\begin{figure}[ht]
\centerline{\includegraphics[width=7cm]{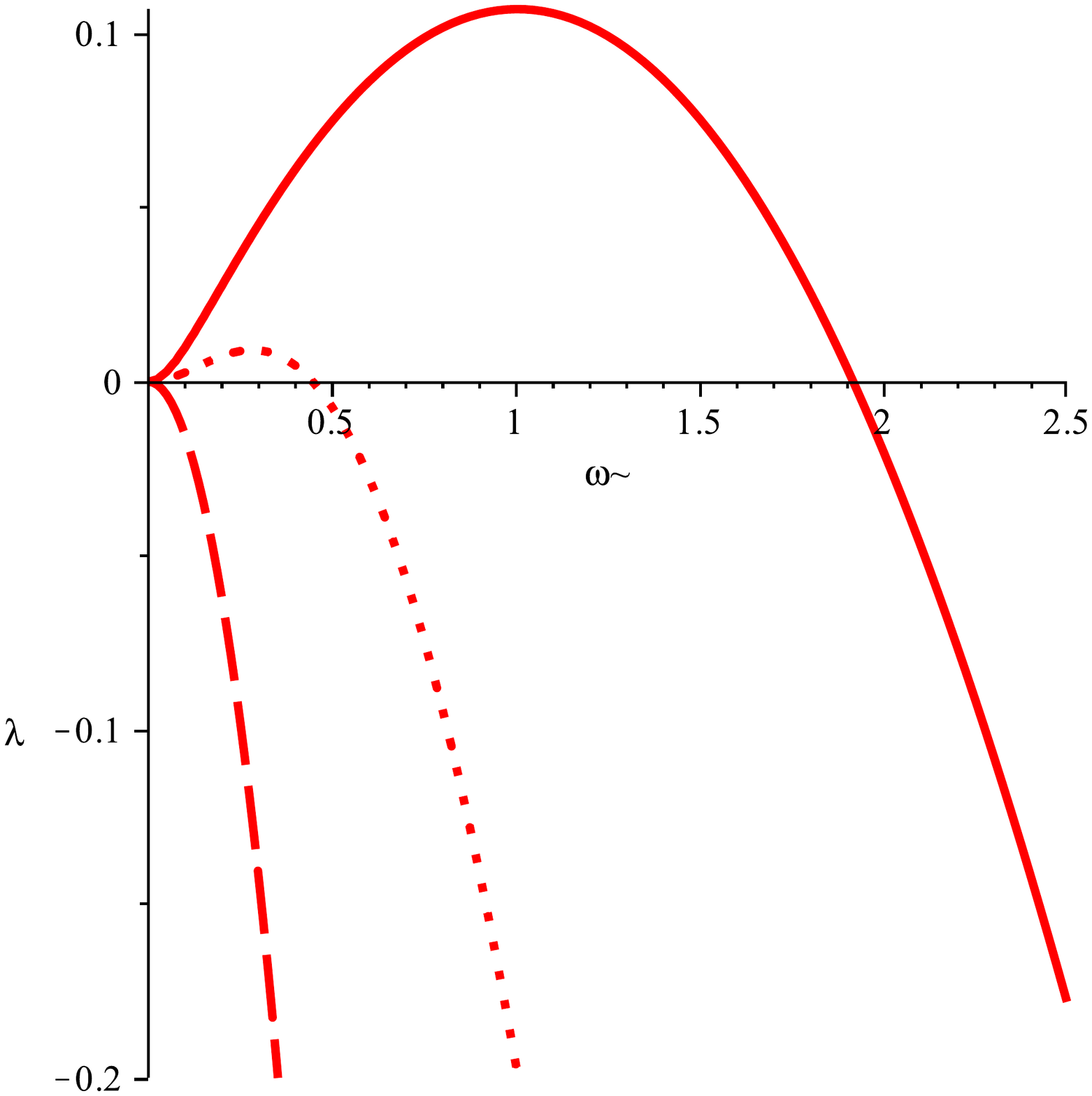}
\hspace{10pt}\includegraphics[width=7cm]{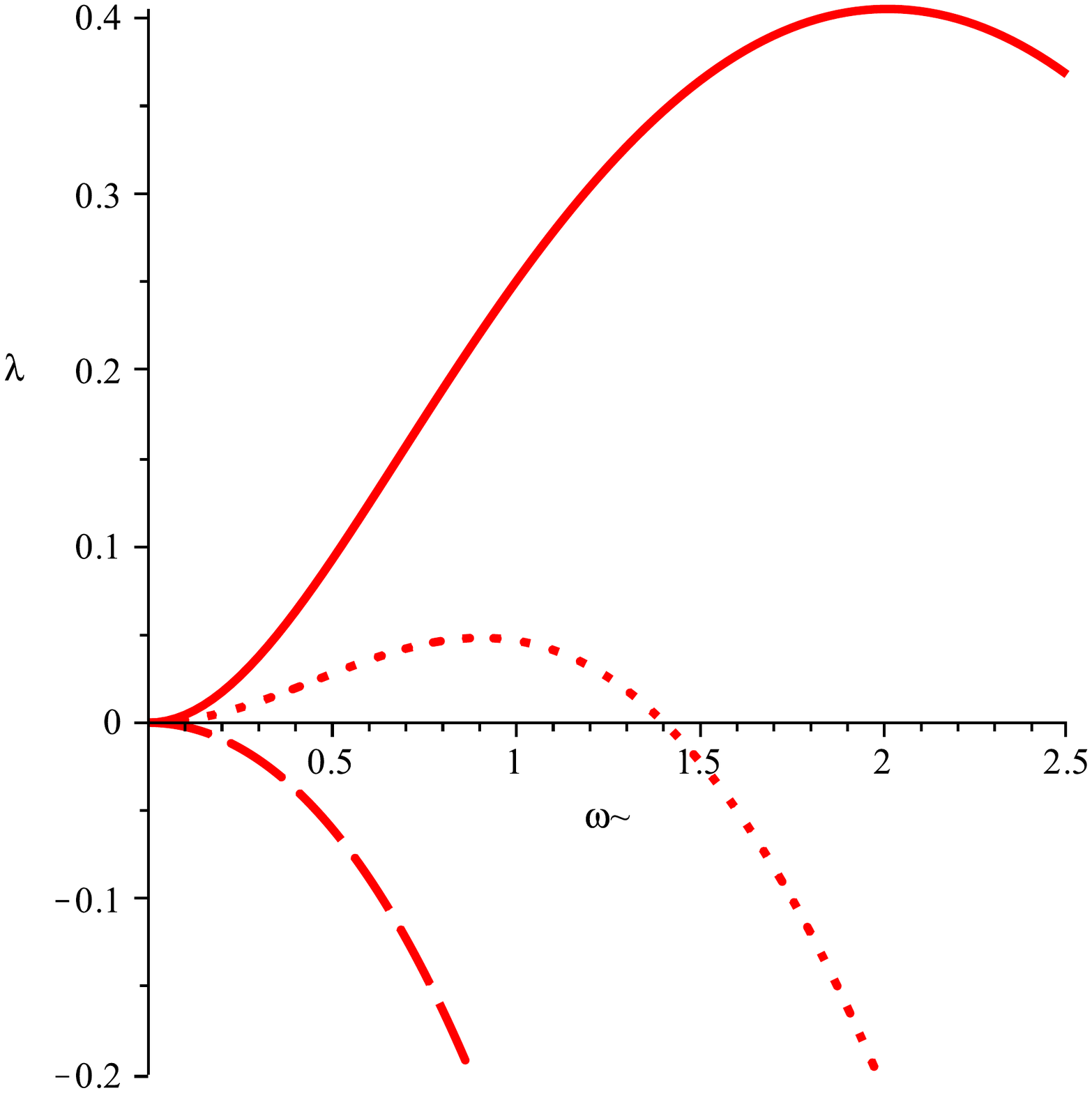}}
\vspace{-12mm}
\caption{\small \sl  (Left) Plot of the dispersion curve $\lambda(\omega)$  for $\theta=0.4$ and $\alpha=0.1$ (continuous  line), $\alpha=0.4$ (dotted line) and $\alpha=1.3$ (dashed line). (Right) Plot of the dispersion curve for $\theta=0.1$ and $\alpha=0.1$ (continuous  line) $\alpha=0.2$ (dotted line) $\alpha=0.4$ (dashed line).}
\label{ever fall in love}
\end{figure}

The intention of this note is to shed some light on the mechanism of the onset of  instability of traveling waves in higher dimension, and in particular to get some idea about the shape of the dispersion curves for more general problems of KPP and ignition type.  We chose to study the planar waves for a simple problem where explicit solutions are available  since, unlike for example in the case of some activator-inhibitor systems (see  \cite{MR1257144,taniguchi_2003,MR1750117}), there does not seem to exist a well established methodology to deal with this issue.  Indeed, the usual approach, involving some limit procedure,  is based on the fact that of one of the components of the system becomes more  concentrated in space, for example it has a form of a spike or undergoes a sharp transition, as the small parameter tends to $0$.  This leads in many cases to a limiting problem for which the spectrum can be completely understood. However, for the KPP or ignition type nonlinearities it is not immediately clear what should the limiting problem be.  We believe that the  instability of the planar fronts described  here is a robust  phenomenon with respect  to change of the nonlinearities. 

\end{remark}

%%%%%%%%%%%%%%%%%%%%%%%%%%%%%%%%%%%
%
%%%%%% BIBLIO %%%%%%%%%%%%%%%%%%%%%%
%
%%%%%%%%%%%%%%%%%%%%%%%%%%%%%%%%%%%%

 \bibliography{BM_unstability.bib}
 \bibliographystyle{plain}
 %%%%%%%%%%%%%%%%%%%%%%%%%%%%%%%%%%%%%%%%%%%%%%%%%%%

%\pagestyle{myheadings}

%\begin{thebibliography}{99}

%\bibitem{BenAmar} Ciarletta, P.;  Foret, L. and Ben Amar, M.  
%\emph{ The radial growth phase of malignant melanoma~: muti-phase modelling, numerical simulation and linear stability}, J. R. Soc. Interface  (2011) {\bf 8} (56), 345--368.

%\bibitem {FreeFem} \href{http:/www.freefem.org/}{http://www.freefem.org} 

%\bibitem{BenJ} Golding, I.; Kozlovsky, Y.; Cohen, I. and Ben~Jacob, E.
%\emph{Studies of bacterial branching growth using reaction--diffusion models for colonial development}, 
% Physica A, 260:510--554, 1998.

%\bibitem{Kessler_Levine} Kessler, D. A. and Levine, H. \emph{Fluctuation-induced diffusive instabilities}; Letters to Nature, Vol.~394 (1998) 556--558.

%\end{thebibliography}

\end{document}